\documentclass{article}
\usepackage{latexsym}
\usepackage{amssymb}
\usepackage{amsmath}
\usepackage{amscd}
\usepackage[only,twlrm,ninrm,sevrm]{rawfonts}
\usepackage[latin1]{inputenc}
\numberwithin{equation}{section} \textwidth=5.45in
\newtheorem{Th}{Theorem}[section]

\newtheorem{Le}{Lemma}[section]

\begin{document}
\title{Existence of least energy nodal solution with two nodal domains for a generalized Kirchhoff problem in an Orlicz Sobolev space}
\author{
Giovany M. Figueiredo\thanks{Partially supported by CNPq/PQ
301242/2011-9 and 200237/2012-8 }  \\
\noindent Universidade Federal do Par\'a, Faculdade de Matem\'atica, \\
\noindent CEP: 66075-110, Bel\'em - Pa, Brazil \\
\noindent e-mail: giovany@ufpa.br  \vspace{0.5cm}\\
\noindent
Jefferson A. Santos \thanks{Partially supported by CNPq-Brazil grant Casadinho/Procad 552.464/2011-2 } \\
\noindent Universidade Federal de Campina Grande,\\
\noindent Unidade Acad\^emica de Matem\'atica e Estat\'istica,\\
\noindent CEP:58109-970, Campina Grande - PB, Brazil\\
\noindent e-mail: jefferson@dme.ufcg.edu.br\vspace{0.5cm}\\
}\vspace{0.5cm}
\date{}

\pretolerance10000

\maketitle

\begin{abstract}
  We show the existence of a nodal solution with two nodal domains for
a generalized Kirchhoff equation of the type
$$
-M\left(\displaystyle\int_\Omega \Phi(|\nabla u|)dx\right)\Delta_\Phi u = f(u) \ \
\mbox{in} \ \ \Omega, \ \ u=0 \ \ \mbox{on} \ \ \partial\Omega,
$$
where $\Omega$ is a bounded domain in $\mathbf{R}^N$, $M$ is a
general $C^{1}$ class function, $f$ is a superlinear $C^{1}$
class function with subcritical growth,  $\Phi$ is defined for $t\in \mathbf{R}$ by setting
$ \Phi(t)=\int_0^{|t|}\phi(s)sds$, $\Delta_\Phi$ is the operator $\Delta_\Phi u:=div(\phi(|\nabla u|)\nabla u)$. The proof is based on 
a minimization argument and a quantitative deformation
lemma.
\end{abstract}

\maketitle


\section{Introduction}

In this paper we study the existence of a nodal solution with two
nodal domains for the following problem
$$
 \left\{
\begin{array}{rcl}
-M\left(\displaystyle\int_\Omega \Phi(|\nabla u|)dx\right)\Delta_\Phi u = f(u) \ \ \mbox{in} \ \ \Omega,\\
u^{+}\neq 0 \ \ \mbox{and} \ \ u^{-}\neq 0  \ \ \mbox{in} \ \ \Omega,\\
u=0 \ \ \mbox{on} \ \partial\Omega,
\end{array}
\right.\leqno{(P)}
$$
where $\Omega \subset \mathbf{R}^{N}$ is a smooth bounded domain and
$$
u^{+}(x)=\max\{u(x),0\} \ \ \mbox{and} \ \ u^{-}(x)=\min\{ u(x),0\},
\ \ \mbox{for all} \ \ x \in \Omega.
$$
Notice that, in this case, $u=u^{+}+u^{-}$ and $|u|=u^{+}-u^{-}$. We
looking for solution with exactly two nodal domains. Problem $(P)$  with $\phi(t)=2$, that is,
$$
\left \{ \begin{array}{l}
-M\left(\displaystyle\int_\Omega \mid\nabla u\mid^{2} dx\right)\Delta u =f(u) \ \mbox{in}\ \Omega,\\
u\in H^{1}_{0}(\Omega)                                    
\end{array}
\right.\leqno{(*)}
$$
is called Kirchhoff type because of the presence of the term
$M\left(\displaystyle\int_{\Omega}|\nabla u|^{2} dx \right)$. Indeed, this operator
appears in the Kirchhoff
equation, which arises in nonlinear vibrations, namely
$$
\left\{
\begin{array}{l}
 u_{tt}-M\left(\displaystyle\int_{\Omega}|\nabla u|^{2} dx \right)\Delta u =
g(x,u) \ \mbox{in} \ \Omega \times (0,T)\\
u=0 \ \mbox{on} \ \partial\Omega \times (0,T)\\
u(x,0)=u_{0}(x)\ \ , \ \ u_{t}(x,0)=u_{1}(x).
\end{array}
\right.
$$

The reader may consult \cite{alvescorrea}, \cite{alvescorreama}, \cite{Cristian}, \cite{ma}
 and the references therein, for more physical motivation on Kirchhoff problem.\\

Before stating our main results, we need the following hypotheses
on the function $M$.

The  function  $M: \mathbf{R}_{+} \to \mathbf{R}_{+}$ is $C^{1}$
class and satisfies the following conditions:

\begin{description}

\item[($M_{1}$)] The function $M$ is increasing and $0 < M(0)=: \sigma_{0}$.

\item[($M_{2}$)] The function $t\mapsto\displaystyle\frac{M(t)}{t}$ is decreasing.

\end{description}

A typical example of function verifying the assumptions
$(M_{1})-(M_{2})$ is given by
$$
\displaystyle M(t)=m_{0}+bt, \ \ \mbox{where} \ \ m_{0}>0  \ \ \mbox{and} \ \
b>0.
$$
This is the example that was considered in \cite{kirchhoff},
\cite{Mao},  \cite{Mao1},  \cite{Shuai} and \cite{Zhang}. More generally, each
function of the form
$$
\displaystyle M(t)=m_{0}+bt+\displaystyle\sum_{i=1}^{k}b_{i}t^{\gamma_{i}}
$$
with $b_{i}\geq 0$ and $\gamma_{i}\in (0,1)$ for all $i\in
\{1,2,\ldots, k\}$ verifies the hypotheses $(M_{1})-(M_{2})$. An
another example is $M(t) = m_{0} + ln(1 + t)$.

\vspace{.3cm}

The hypotheses on the function $\phi:(0,+\infty)\rightarrow
(0,+\infty)$ of $C^{1}$ class are the following:

\begin{description}

\item[($\phi_1$)] For all $t>0$, 
$$
\phi(t)>0 \ \ \mbox{and} \ \ (\phi(t)t)'>0.
$$

\item[($\phi_2$)] There exist $l \in (\frac{N}{2},N)$,  $l< m < \min\{\frac{l^{*}}{2}, N\}$, where $l^{*}=
\displaystyle\frac{lN}{(N-l)}$ such that
$$
l\leq \frac{\phi(t)t^{2}}{\Phi(t)}\leq m,
$$
for $t> 0$, where $\Phi(t)=\displaystyle\int^{|t|}_{0}\phi(s)s ds$.

\item[($\phi_3$)] For all $t>0$, 
$$l-1\leq  \frac{(\phi(t)t)'t}{\phi(t)t}\leq m-1.$$
\end{description}

We assume that the function 
$f$ is $C^{1}$ class and satisfies
\begin{description}
\item[($f_{1}$)] $$\lim_{|t|\to 0^{+}}\frac{f(t)}{\phi(t)t}=0 $$
and
\item[($f_{2}$)] 
$$
\lim_{|t|\to \infty}\frac{f(t)}{\phi_{*}(t)t}=0.
$$

\item[($f_{3}$)]  There is $\theta\in (2m,l^{*})$ such that
$$
0<\theta F(t)\leq f(t)t, \ \forall |t|>0, \ \ \mbox{where} \ \
F(t)=\displaystyle\int_{0}^{t}f(s)ds.
$$

\item[($f_{4}$)] The map $$t\mapsto\frac{f(t)}{t^{2m-1}}$$ is increasing in $|t|>0$.

\end{description}

The main result of this paper is:

\begin{Th}\label{Theorem1}
Suppose that $(M_{1})-(M_{2})$, $(\phi_{1})-(\phi_{3})$  $(f_{1})-(f_{4})$ hold. Then problem $(P)$
possesses a least energy nodal solution, which has precisely two
nodal domains.
\end{Th}

In the last years, many authors show the existence and/or multiplicity of nontrivial solutions for the 
problem $(*)$, as can be seen \cite{AlvesFigueiredo}, \cite{Azzollini}, \cite{Cammaroto}, \cite{Cheng}, \cite{Chen}, \cite{SChen},
\cite{jmaa}, \cite{Cristian}, \cite{Giovany1}, \cite{Giovany2}, \cite{He}, \cite{Li}, \cite{Liao}, \cite{Liu}, \cite{Liu1}, \cite{Naimen}, 
 \cite{Xiang}, \cite{Wang1}, \cite{Wang2} and reference therein.\\

Only the articles \cite{Giovany3}, \cite{Mao}, \cite{Mao1}, \cite{Shuai} and \cite{Zhang} consider solutions that change
sign (nodal solution) for the Kirchhoff problem. In \cite{Mao}, \cite{Mao1}, \cite{Shuai} and \cite{Zhang} the authors show the
existence of solutions which change sign considering $M(t)=a+bt$ for some positive constants $a,b$. In \cite{Zhang} 
the authors showed the existence of the sign
changing solutions to the problem $(*)$ for the cases
$4$-sublinear, asymptotically $4$-linear and $4$-superlinear. In
that paper the authors use variational methods and invariant sets of
descent flow. In \cite{Mao} the authors showed the same result found
in \cite{Zhang} without considering the Ambrosetti-Rabinowitz condition. In \cite{Mao1} the authors study the result found in
\cite{Zhang} considering now the case asymptotically $3$-linear.  In \cite{Giovany3} and \cite{Shuai} the authors show the existence of nodal solution
using a minimization argument and a quantitative deformation lemma. In \cite{Giovany3} the function $M$ was more general 
that included the case $M(t)=a+bt$.\\
 
The generalized  Kirchhoff problem, that is, Kirchhoff problem in Orlicz-Sobolev spaces, was studied in \cite{Chung}, 
\cite{Giovany4} and \cite{Giovany5} using Genus theory in order to show a multiplicity result. In \cite{Chung1} 
the author used the Mountain Pass  Theorem to show the existence result.  In \cite{Chung2} the author used the 
Ricceri's three points critical result to show multiplicity result.\\

In this work we completes the results found in \cite{Chung}, \cite{Chung1}, \cite{Chung2},  \cite{Giovany4}, \cite{Giovany5}  and extend the studies found in  
\cite{Giovany3},  \cite{Mao}, \cite{Mao1}, \cite{Shuai} and \cite{Zhang} in the
following sense:

\noindent a)  Unlike \cite{Chung}, \cite{Chung1}, \cite{Chung2}, \cite{Giovany4}, \cite{Giovany5}, we show the existence of a nodal solution 
for a generalized  Kirchhoff problem.  \\

\noindent b) Since we work in Orlicz-Sobolev spaces, some estimates were necessary, as for example, Lemma 3.4 and some results were more 
delicate, as for example, Lemma 3.5. \\

\noindent c) This result is new even in the case $M \equiv 1$. \\

\noindent d) Problem $(P_\alpha)$ possesses more complicated
nonlinearities, as for example: 

\noindent (i) $\Phi(t)=t^{p_{0}}+ t^{p_{1}}$, $1< p_{0} < p_{1}< N$
and $ p_1\in (p_0,p^{*}_{0})$, for $N \geq 2$.

\noindent (ii) $\Phi(t)=(1+t^{2})^{\gamma}-1$ for $\gamma\in (1,3)$ and $N=3$.

\noindent (iii) $\Phi(t)=t^p log(t^q+1)$, with $p\in (\frac{N}{2},N)$, $q\in (0,\min\{\frac{p}{N-p},N-p\})$ and 
$N\geq 3$.\\

The paper is organized as follows. In the next section we give a brief review on Orlicz-Sobolev spaces. 
In section 3 we give the variational framework and we prove some technical lemmas. 
In the section 4 we prove Theorem \ref{Theorem1}.


\section{A brief review on Orlicz-Sobolev spaces}\label{Section Orlicz}

Let $\varphi$ be a real-valued function defined in $[0,\infty)$ and
having the following properties:

\noindent $a)$ \ \  $\varphi(0)=0$, $\varphi(t)>0$ if $t>0$ and
$\displaystyle\lim_{t\rightarrow \infty}\varphi(t)=\infty$.

\noindent $b)$ \ \  $\varphi$ is nondecreasing, that is, $s>t$ implies
$\varphi(s) \geq \varphi(t)$.

\noindent $c)$ \ \ $\varphi$ is right continuous, that is,
$\displaystyle\lim_{s\rightarrow t^{+}}\varphi(s)=\varphi(t)$.

Then, the real-valued function $\Phi$ defined on $\mathbf{R}$ by
$$
\Phi(t)= \displaystyle\int^{|t|}_{0}\varphi(s) \ ds
$$
is called an N-function. For an N-function $\Phi$ and  an open set
$\Omega \subseteq \mathbf{R}^{N}$, the Orlicz space
$L_{\Phi}(\Omega)$ is defined (see \cite{adams} or \cite{Gossez}). When $\Phi$
satisfies $\Delta_{2}$-condition, that is, when there are $t_{0}\geq
0$ and $K>0$ such that $\Phi(2t)\leq K\Phi(t)$, for all $t\geq
t_{0}$,  the space $L_{\Phi}(\Omega)$ is the vectorial space of the
measurable functions $u: \Omega \to \mathbf{R}$ such that
$$
\displaystyle\int_{\Omega}\Phi(|u|) \ dx < \infty.
$$
The space $L_{\Phi}(\Omega)$ endowed with Luxemburg norm, that is,
the norm given by
$$
|u|= \inf \biggl\{\tau >0:
\int_{\Omega}\Phi\Big(\frac{|u|}{\tau}\Big)\ dx\leq 1\biggl\},
$$
is a Banach space. The complement function of $\Phi$, denoted by
$\widetilde{\Phi}$, is given by the Legendre transformation, that is
$$
\widetilde{\Phi}(s)=\displaystyle\max_{t \geq 0}\{st -\Phi(t)\} \ \
\mbox{for} \ \ s \geq 0.
$$
These $\Phi$ and $\widetilde{\Phi}$ are complementary each other.
Involving the functions $\Phi$ and $\widetilde{\Phi}$, we have the
Young's inequality given by
$$
st \leq \Phi(t) + \widetilde{\Phi}(s).
$$
Using the above inequality, it is possible to prove the following
H\"{o}lder type inequality
$$
\biggl|\displaystyle\int_{\Omega}u v \ dx \biggl|\leq
2|u|_{\Phi}|v|_{\widetilde{\Phi}}\,\,\, \forall \,\, u \in
L_{\Phi}(\Omega) \,\,\, \mbox{and} \,\,\,  v \in
L_{\widetilde{\Phi}}(\Omega).
$$

Hereafter, we denote by $W^{\Phi}_{0}(\Omega)$ the
Orlicz-Sobolev space obtained by the completion of
$C^{\infty}_{0}(\Omega)$ with norm
$$
\|u\|=|u|_\Phi+|\nabla u|_\Phi.
$$

When $\Omega$ is bounded, for all $u \in  W^{\Phi}_{0}(\Omega)$, there is $c>0$ such that
\begin{eqnarray}\label{Poincare}
|u| \leq c |\nabla u|_\Phi
\end{eqnarray}
and 
\begin{eqnarray}\label{Poincare1}
\displaystyle\int_{\Omega}\Phi(u) dx\leq c \displaystyle\int_{\Omega}\Phi(|\nabla u|)dx.
\end{eqnarray}
In this case, we can consider

$$
\|u\|=|\nabla u|_{\Phi}.
$$

Another important function related to function $\Phi$, is the
Sobolev conjugate function $\Phi_{*}$ of $\Phi$ defined by
$$
\Phi^{-1}_{*}(t)=\displaystyle\int^{t}_{0}\displaystyle\frac{\Phi^{-1}(s)}{s^{(N+1)/N}}ds, \ t>0.
$$

Another important inequality  was proved by
Donaldson and Trudinger \cite{Donaldson2}, which establishes that
for all open $\Omega \subset \mathbf{R}^{N}$ and there is a constant  $S_N=S(N) > 0$ such that
\begin{equation}\label{trudinger-emb}
|u|_{\Phi_*}\leq S_N|\nabla u|_{\Phi}, \ \ \mbox{for all} \ \ u\in W_0^{1,\Phi}(\Omega).
\end{equation}
This inequality shows the below embedding is continuous
$$
W_0^{1,\Phi}(\Omega) \stackrel{\hookrightarrow}{\mbox{\tiny cont}} L_{\Phi_*}(\Omega).
$$
If bounded domain $\Omega$ and the limits below hold
\begin{equation} \label{M1}
\limsup_{t \to 0}\frac{B(t)}{\Phi(t)}< +\infty \,\,\, \mbox{and} \,\,\, \limsup_{|t| \to +\infty}\frac{B(t)}{\Phi_{*}(t)}=0,
\end{equation}
the embedding
\begin{equation} \label{M2}
W_0^{1,\Phi}(\Omega) \hookrightarrow L_{B}(\Omega)
\end{equation}
is compact.

The hypotheses $(\phi_{1})-(\phi_{3})$ implies that $\Phi$,
$\widetilde{\Phi}$, $\Phi_{*}$ and $\widetilde{\Phi}_{*}$ satisfy
$\Delta_{2}$-condition. This condition allows us conclude that:

\noindent 1) $u_{n}\rightarrow 0$ in $L_{\Phi}(\Omega)$ if, and only
if, $\displaystyle\int_{\Omega}\Phi(u_{n})\ dx\rightarrow 0$.

\noindent 2)  $L_{\Phi}(\Omega)$ is separable and
$\overline{C^{\infty}_{0}(\Omega)}^{|.|_{\Phi}}=L_{\Phi}(\Omega)$.

\noindent 3)  $L_{\Phi}(\Omega)$ is reflexive and its dual is
$L_{\widetilde{\Phi}}(\Omega)$(see \cite{adams}).

\vspace{.5cm}

Under assumptions $(\phi_{1})-(\phi_{3})$, some elementary
inequalities listed in the following lemmas are valid. For the
proofs, see \cite{fukagai}.

\begin{Le}
Assume $(\phi_1)$ and $(\phi_2)$. Then,
$$
\Phi(t) = \int_0^{|t|} s \phi(s) ds, \,\,\,
$$
is a $N$-function with $\Phi, \widetilde{\Phi} \in \Delta_2$.  Hence, $L_\Phi(\Omega)$ and $W_0^{1,\Phi}(\Omega)$ are reflexive and separable spaces.
\end{Le}

\begin{Le}\label{DESIGUALD}The functions $\Phi$, $\Phi_*$,  $\widetilde{\Phi}$ and $\widetilde{\Phi}_*$ satisfy the inequalities
\begin{equation} \label{D1}
\widetilde{\Phi}(\phi(t)t) \leq \Phi(2t)  \mbox{ and } \widetilde{\Phi}\bigl(\frac{\Phi(t)}{t}\bigl) \leq \Phi(t)\,\,\, \forall t \geq 0.
\end{equation}
\end{Le}

\begin{Le} \label{desigualdadeimportantes} Assume $(\phi_1)$ and $(\phi_2)$ hold and let $\xi_{0}(t)=\min\{t^{l},t^{m}\}$,\linebreak $ \xi_{1}(t)=\max\{t^{l},t^{m}\},$ for all $t\geq 0$. Then,
$$
\xi_{0}(\rho)\Phi(t) \leq \Phi(\rho t) \leq  \xi_{1}(\rho)\Phi(t) \;\;\; \mbox{for} \;\; \rho, t \geq 0
$$
and
$$
\xi_{0}(|u|_{\Phi}) \leq \int_{\Omega}\Phi(u) dx\leq \xi_{1}(|u|_{\Phi})  \;\;\; \mbox{for} \;\; u \in L_{\Phi}(\Omega).
$$
\end{Le}

\begin{Le} \label{novomasseguefukagai} Assume $(\phi_3)$ holds and let $\Psi(t)=\phi(t)t$, $\tau_{0}(t)=\min\{t^{l-1},t^{m-1}\}$,
\linebreak $ \tau_{1}(t)=\max\{t^{l-1},t^{m-1}\},$ for all $t\geq 0$. Then,
$$
\tau_{0}(\rho)\Psi(t) \leq \Psi(\rho t) \leq  \tau_{1}(\rho)\Psi(t) \;\;\; \mbox{for} \;\; \rho, t \geq 0.
$$
\end{Le}

\begin{Le} \label{F3} The function $\Phi_*$ satisfies the following inequality
\begin{eqnarray}\label{desigualdadecritica}
l^{*} \leq \frac{\Phi'_*(t)t}{\Phi_{*}(t)} \leq m^{*} \,\,\, \mbox{for} \,\,\, t > 0.
\end{eqnarray}
\end{Le}
As an immediate consequence of the Lemma \ref{F3}, we have the following result:

\begin{Le} \label{F2} Assume $(\phi_1)-(\phi_2)$ hold and let  $\xi_{2}(t)=\min\{t^{l^{*}},t^{m^{*}}\},$ $\xi_{3}(t)=\max\{t^{l^{*}},t^{m^{*}}\}$ for all $t\geq 0$. Then,
$$
\xi_{2}(\rho)\Phi_*(t) \leq \Phi_*(\rho t) \leq \xi_{3}(\rho)\Phi_*(t) \;\;\; \mbox{for} \;\; \rho, t \geq 0
$$
and
$$
\xi_{2}(|u|_{\Phi_*}) \leq \int_{\Omega}\Phi_*(u)dx \leq \xi_{3}(|u|_{\Phi_*})  \;\;\; \mbox{for} \;\; u \in L_{\Phi_*}(\Omega).
$$
\end{Le}
\begin{Le} \label{lem Phiest}
Let $\widetilde{\Phi}$ be the complement of $\Phi$ and put
$$
\xi_4(s)=\min\{s^{\frac{l}{l-1}}, s^{\frac{m}{m-1}}\}\ \mbox{and}\
\xi_5(s)=\max\{s^{\frac{l}{l-1}}, s^{\frac{m}{m-1}}\}, \ s\geq0.
$$
Then the following inequalities hold
$$
\xi_4(r)\widetilde{\Phi}(s)\leq\widetilde{\Phi}(r s)\leq
\xi_5(r)\widetilde{\Phi}(s),\ r,s\geq0
$$
and
$$
\xi_4(| u|_{\widetilde{\Phi}})\leq
\int_\Omega\widetilde{\Phi}(u)dx\leq \xi_5(|
u|_{\widetilde{\Phi}}),\ u\in
L_{\widetilde{\Phi}}(\Omega).
$$
\end{Le} \vskip0.5cm


\section{Variational framework and technical lemmas}

We say that $u \in W^{1,\Phi}_{0}(\Omega)$ is a weak nodal solution of
the problem $(P)$ if $u^{+}\neq 0$, $u^{-}\neq 0$ in $\Omega$ and it
verifies
$$
M\biggl(\displaystyle\int_{\Omega}\Phi(|\nabla u|\ dx \biggl)\displaystyle\int_{\Omega} \phi(|\nabla u|)\nabla u\nabla v \ dx -
\displaystyle\int_{\Omega}f(u)v \ dx=0, \ \ \mbox{for all} \ \
v \in W^{1,\Phi}_{0}(\Omega).
$$

In view of $(f_1)-(f_2)$, we have that the functional
$J:W^{1,\Phi}_{0}(\Omega)\to \mathbf{R}$ given by
$$
J(u) := \widehat{M}\biggl(\displaystyle\int_{\Omega}\Phi(|\nabla u|\ dx \biggl)
 -\displaystyle\int_{\Omega} F(u) \ dx
$$
is well defined, where $\widehat{M}(t)=\displaystyle\int^{t}_{0}M(s)
\ ds$. Moreover, $J \in C^1(W^{1,\Phi}_{0}(\Omega),\mathbf{R})$ with the
following derivative
$$
J'(u)v=M\biggl(\displaystyle\int_{\Omega}\Phi(|\nabla u|)\ dx \biggl)\displaystyle\int_{\Omega} \phi(|\nabla u|)\nabla u\nabla v \ dx -
\displaystyle\int_{\Omega} f(u)v \ dx, 
$$
for all $v \in W^{1,\Phi}_{0}(\Omega)$. Thus, the weak solutions of $(P)$ are  precisely the critical points
of $J$. Associated to the functional $J$ we define the Nehari
manifold
$$
\mathcal{N}:= \biggl\{u\in W^{1,\Phi}_{0}(\Omega)\backslash\{0\}:J'(u)u=0
\biggl\}.
$$

In the Theorem \ref{Theorem1} we prove that there is $w \in
\mathcal{M}$ such that
$$
J(w)= \displaystyle\min_{v \in \mathcal{M}}J(v),
$$
where
$$
\mathcal{M}:= \biggl\{w\in \mathcal{N}: J'(w)w^{+}=0= J'(w)w^{-}
\biggl\}.
$$

From $(M_{1})$ we have that $M\biggl(\displaystyle\int_{\Omega}\Phi(|\nabla w^{\pm}|)\ dx \biggl)\leq 
M\biggl(\displaystyle\int_{\Omega}\Phi(|\nabla w|) \ dx \biggl)$,
for $w \in \mathcal{M}$. Thus, this last inequality implies that
\begin{eqnarray}\label{beleza}
J'(w^{\pm})w^{\pm}\leq 0, \ \ \mbox{for all} \ \ w \in \mathcal{M}.
\end{eqnarray}

Let us begin by establishing some preliminary results which will be
exploited  in the last section for a minimization argument.

\begin{Le} \label{limitacaoporbaixo1}
\begin{description}
\item[(a)] For all $ u \in \mathcal{N}$ we have
$$
J(u)\geq \frac{(\theta -4m)}{4\theta}\sigma_{0}\xi_{0}(|\nabla u|_{\Phi}).
$$

\item[(b)]  There is $\rho>0$ such that
$$
\|u\|\geq \rho, \ \ \mbox{for all} \ \ u \in \mathcal{N}
$$
and
$$
\|w^{\pm}\|\geq \rho, \ \ \mbox{for all} \ \ w \in \mathcal{M}.
$$
\end{description}
\end{Le}
\noindent {\bf Proof.} From definition of $\widehat{M}$ and
$(M_{2})$, we get
\begin{eqnarray}\label{estimativaM1}
\widehat{M}(t)\geq \frac{1}{2}M(t)t, \ \ \mbox{for all}\ \ t\geq 0.
\end{eqnarray}

Thus, by $(f_{3})$ we get
$$
J(u)= J(u)-\frac{1}{\theta}J'(u)u\geq \widehat{M}\biggl(\displaystyle\int_{\Omega}\Phi(|\nabla u|)dx\biggl)
-\frac{1}{\theta}M\biggl(\displaystyle\int_{\Omega}\Phi(|\nabla u|)dx\biggl)\displaystyle\int_{\Omega}\phi(|\nabla u|)|\nabla u|^{2} dx.
$$

Using $(\phi_{2})$ we obtain
$$
J(u)\geq \widehat{M}\biggl(\displaystyle\int_{\Omega}\Phi(|\nabla u|)dx\biggl)
-\frac{m}{\theta}M\biggl(\displaystyle\int_{\Omega}\Phi(|\nabla u|)dx\biggl)\displaystyle\int_{\Omega}\Phi(|\nabla u|)dx.
$$

Since $\theta \in (2m,l^{*})$, by (\ref{estimativaM1}) and Lemma \ref{desigualdadeimportantes}, the proof of $(a)$ is finished.

To prove $(b)$, notice that  for all $u \in W^{1,\Phi}_{0}(\Omega)$,  by $(f_{1})$, $(f_{2})$, $(\phi_{2})$, (\ref{Poincare1}) 
and (\ref{desigualdadecritica}), given
$\epsilon>0$, there exists $C_{\epsilon}>0$ such that
\begin{eqnarray}\label{estimativaf1}
\displaystyle\int_{\Omega}f(u)u dx  \leq \epsilon c m \displaystyle\int_{\Omega}\Phi(|\nabla u|)  dx
+ C_{\epsilon}cm^{*}\displaystyle\int_{\Omega}\Phi_{*}( u)dx.
\end{eqnarray}

Now from definition of $\mathcal{N}$, $(M_{1})$, $(\phi_{2})$ again and (\ref{estimativaf1}),   we have
$$
(\sigma_{0}l - \epsilon c m)\displaystyle\int_{\Omega}\Phi(|\nabla u|)  dx 
\leq C_{\epsilon} cm^{*}\displaystyle\int_{\Omega}\Phi_{*}( u)dx.
$$

By Lemma \ref{desigualdadeimportantes} and (\ref{trudinger-emb}) we obtain 
$$
(\sigma_{0}l - \epsilon cm)\xi_{0}(|\nabla u|_{\Phi})  
\leq cCC_{\epsilon}m \xi_{3}(|\nabla u|_{\Phi}).
$$
If $|\nabla u|_{\Phi}\geq 1$, the proof is done. Suppose that $|\nabla u|_{\Phi}\leq 1$. 
Then, from the last inequality we have
$$
(\sigma_{0}l - \epsilon cm)|\nabla u|_{\Phi}^{m}  
\leq cCC_{\epsilon}m |\nabla u|_{\Phi}^{l^{*}}.
$$

Consequently, 
$$
\biggl(\frac{\sigma_{0}l - \epsilon c m}{cCC_{\epsilon}m}\biggl)^{\frac{1}{(l^{*}-m)}}
\leq|\nabla u|_{\Phi},
$$
for all $u \in \mathcal{N}$. Then proof of $(b)$ is finished with 
$\rho=\min\biggl\{1, \biggl(\frac{\sigma_{0}l - \epsilon c m}
{cCC_{\epsilon}m}\biggl)^{\frac{1}{(l^{*}-m)}}\biggl\}$.

From (\ref{beleza}) and repeating the reasoning before we obtain
$$
0< \rho\leq \|w^{\pm}\|. \hfill\rule{2mm}{2mm}
$$

We apply the next result in the last section   to every bounded
minimizing sequence of $J$ on $\mathcal{M}$ in order to ensure that
the candidate minimizer is different from zero.

\begin{Le} \label{limitacaoporbaixo2}
If $(w_{n})$ is a bounded sequence in $\mathcal{M}$, then there is a positive constant $q \in (m, l^{*})$ such that
$$
\displaystyle\liminf_{n\rightarrow
\infty}\displaystyle\int_{\Omega}|w_{n}^{\pm}|^{q} \ dx >0.
$$
\end{Le}
\noindent {\bf Proof.} Note that  by $(f_{1})$, $(f_{2})$, $(\phi_{2})$ 
and (\ref{desigualdadecritica}), given
$\epsilon>0$, there exists $C_{\epsilon}>0$ such that
\begin{eqnarray}\label{estimativaf11}
\displaystyle\int_{\Omega}f(w_{n}^{\pm})w_{n}^{\pm} dx  \leq \epsilon m \displaystyle\int_{\Omega}\Phi(|\nabla w_{n}^{\pm}|)  dx
+ \epsilon m^{*}\displaystyle\int_{\Omega}\Phi_{*}( w_{n}^{\pm})dx + C_{\epsilon}\displaystyle\int_{\Omega}|w_{n}^{\pm}|^{q} \ dx.
\end{eqnarray}

On the other hand, by $(M_{1})$, $(\phi_{2})$ and Lemma \ref{desigualdadeimportantes}, we get
\begin{eqnarray}\label{estimativaf12}
\sigma_{0}l\xi_{0}(|\nabla w_{n}^{\pm}|_{\Phi}) \leq l\sigma_{0}\displaystyle\int_{\Omega}\Phi( |\nabla w_{n}^{\pm}|)dx 
\leq \displaystyle\int_{\Omega}f(w_{n}^{\pm})w_{n}^{\pm} dx.
\end{eqnarray}

Using (\ref{estimativaf11}), (\ref{estimativaf12})  and last Lemma, we obtain
$$
0<\sigma_{0}\xi_{0}(\rho)\leq \epsilon m \displaystyle\int_{\Omega}\Phi(|\nabla w_{n}^{\pm}|)  dx
+ \epsilon m^{*}\displaystyle\int_{\Omega}\Phi_{*}( w_{n}^{\pm})dx + C_{\epsilon}\displaystyle\int_{\Omega}|w_{n}^{\pm}|^{q} \ dx.
$$

Since $(w_{n})$ is bounded in $W^{1,\Phi}_{0}(\Omega)$, by Lemma \ref{desigualdadeimportantes},  there is $C>0$ such that
$$
0<\sigma_{0}l\xi_{0}(\rho)\leq \epsilon C+ C_{\epsilon}
\displaystyle\int_{\Omega}|w_{n}^{\pm}|^{q} \ dx
$$
and the result follows of the last inequality. \hfill\rule{2mm}{2mm}

Next results try to infer geometrical information of $J$ with
respect to $\mathcal{M}$ in the same way that one is used to do
about $\mathcal{N}$. To be more precise, note the similarity between
the next result and that which states that for each $v \in
W^{1,\Phi}_{0}(\Omega) \backslash\{0\}$ there exists $t_{v} > 0$ such that
$t_{v}v \in \mathcal{N}$.

\begin{Le} \label{teovalormedio}
If $v\in W^{1,\Phi}_{0}(\Omega)$ with $v^{\pm}\neq 0$, then there are
$t,s>0$ such that
$$
J'(tv^{+}+sv^{-})v^{+}=0
$$
and
$$
J'(tv^{+}+sv^{-})v^{-}=0.
$$
\end{Le}
\noindent {\bf Proof.} Let
$V:(0,+\infty)\times(0,+\infty)\rightarrow \mathbf{R}^{2}$ be a
continuous function given by
$$
V(t,s)=(J'(tv^{+}+sv^{-})(tv^{+}),J'(tv^{+}+sv^{-})(sv^{-})).
$$
Note that
\begin{eqnarray}\label{prausar1}
J'(tv^{+}+sv^{-})(tv^{+})&=&  M\biggl(\displaystyle\int_{\Omega}\Phi\bigl(|\nabla (tv^{+}+sv^{-})|\bigl)dx\biggl)\!\!\!
\displaystyle\int_{\Omega}\phi \bigl(|\nabla (tv^{+})|\bigl)|\nabla (tv^{+})|^{2} dx \nonumber\\ 
&&-\displaystyle\int_{\Omega}f(tv^{+})tv^{+}dx.
\end{eqnarray}

Using $(M_{1})$ and $(\phi_{2})$ we have
\begin{eqnarray*}
J'(tv^{+}+sv^{-})(tv^{+})\geq  \sigma_{0}l \displaystyle\int_{\Omega}\Phi\bigl(|\nabla (tv^{+})|\bigl)dx
-\displaystyle\int_{\Omega}f(tv^{+})tv^{+}
dx.
\end{eqnarray*}

Considering (\ref{estimativaf1}) in the last inequality we get
\begin{eqnarray*}
J'(tv^{+}+sv^{-})(tv^{+})\geq  \sigma_{0}l \displaystyle\int_{\Omega}\Phi\bigl(|\nabla (tv^{+})|\bigl)dx
-\epsilon\displaystyle\int_{\Omega}\Phi(tv^{+})dx - C_{\epsilon}\displaystyle\int_{\Omega}\Phi_{*}(tv^{+})
dx.
\end{eqnarray*}

From (\ref{Poincare1}) we obtain

\begin{eqnarray*}
J'(tv^{+}+sv^{-})(tv^{+})\geq  (\sigma_{0}l - \epsilon c)\displaystyle\int_{\Omega}\Phi\bigl(|\nabla (tv^{+})|\bigl)dx
- C_{\epsilon}\displaystyle\int_{\Omega}\Phi_{*}(tv^{+})
dx.
\end{eqnarray*}

Now, by Lemma \ref{desigualdadeimportantes} we derive
\begin{eqnarray*}
J'(tv^{+}+sv^{-})(tv^{+})\geq  (\sigma_{0}l - \epsilon c)\xi_{0}(t)\xi_{0}\bigl(|\nabla (v^{+})|_{\Phi}\bigl)
- C_{\epsilon}\xi_{3}(t)\xi_{3}\bigl(|\nabla (v^{+})|_{\Phi}\bigl).
\end{eqnarray*}

Thus, there exists $r>0$ sufficiently small such
that
$$
J'(rv^{+}+sv^{-})(rv^{+})\geq (\sigma_{0}l - \epsilon c)r^{m}\xi_{0}\bigl(|\nabla (v^{+})|_{\Phi}\bigl)-C_{\epsilon}
r^{l^{*}}\xi_{3}\bigl(|\nabla (v^{+})|_{\Phi}\bigl)  >0, \ \ \mbox{for all} \ \ s>0.
$$

Arguing of the same way we get
$$
J'(tv^{+}+rv^{-})(rv^{-})>0, 
$$
for all $t>0$ and  $r>0$ sufficiently small.

On the other hand, using $(\phi_{2})$ in  (\ref{prausar1}) we get
\begin{eqnarray}\label{Gaetano}
J'(tv^{+}+sv^{-})(tv^{+})\leq && m M\biggl(\displaystyle\int_{\Omega}\Phi\bigl(|\nabla (tv^{+})|\bigl)dx\biggl)\!\!
\displaystyle\int_{\Omega}\Phi\bigl(|\nabla (tv^{+}+sv^{-})|\bigl)dx \nonumber\\
&&-\displaystyle\int_{\Omega}f(tv^{+})tv^{+}dx.
\end{eqnarray}

Note that, by $(M_{2})$, there exists $K_{1}>0$ such that
\begin{eqnarray}\label{estimativaporcima}
M(t)\leq M(1)t+K_{1}, \ \ \mbox{for all} \ \ t \geq 0 .
\end{eqnarray}

Using (\ref{estimativaporcima}) in (\ref{Gaetano}) and recalling that $v^{+}$ and $v^{-}$ have compact support disjoint we obtain

\begin{multline*}
J'(tv^{+}+sv^{-})(tv^{+})\leq m M(1)\biggl(\displaystyle\int_{\Omega}\Phi\bigl(|\nabla (tv^{+})|\bigl)dx\biggl)^{2}\\
+mM(1) \biggl(\displaystyle\int_{\Omega}\Phi\bigl(|\nabla (tv^{+})|\bigl)dx\biggl)
\biggl(\displaystyle\int_{\Omega}\Phi\bigl(|\nabla (sv^{-})|\bigl)dx\biggl) \\
+ mK_{1} \displaystyle\int_{\Omega}\Phi\bigl(|\nabla (tv^{+})|\bigl)dx - \displaystyle\int_{\Omega}f(tv^{+})tv^{+}dx.
\end{multline*}

By Lemma \ref{desigualdadeimportantes} and $(f_{3})$ we have 

\begin{multline}\label{Gaetano1}
J'(tv^{+}+sv^{-})(tv^{+})\leq m M(1)\xi_{1}(t)^{2}\xi_{1}(|\nabla (v^{+})|_{\Phi})^{2}\\
+mM(1) \xi_{1}(t)\xi_{1}(s) \xi_{1}(|\nabla (tv^{+})|_{\Phi})
\xi_{1}(|\nabla (v^{-})|_{\Phi}) \\
+ mK_{1} \xi_{1}(t)\xi_{1}(|\nabla (v^{+})|_{\Phi}) - \theta\displaystyle\int_{\Omega}F(tv^{+}) dx.
\end{multline}

Note that, by  $(f_{3})$ again, there are $K_{2}, K_{3}>0$ such that
\begin{eqnarray}\label{estimativaf2}
F(t)\geq K_{2}t^{\theta}-K_{3}.
\end{eqnarray}

Using  (\ref{estimativaf2}) in (\ref{Gaetano1}), we have

\begin{multline*}
J'(tv^{+}+sv^{-})(tv^{+})\leq m M(1)\xi_{1}(t)^{2}\xi_{1}(|\nabla (v^{+})|_{\Phi})^{2}\\
+mM(1) \xi_{1}(t)\xi_{1}(s) \xi_{1}(|\nabla (tv^{+})|_{\Phi})
\xi_{1}(|\nabla (v^{-})|_{\Phi}) \\
+ mK_{1} \xi_{1}(t)\xi_{1}(|\nabla (v^{+})|_{\Phi}) -\theta t^{\theta}K_{2}|v^{+}|^{\theta}_{\theta}+ K_{3}\theta|\Omega|,
\end{multline*}
where $|\Omega|$ denotes the Lebesgue measure of $\Omega$. Thus,
since $\theta>2m$, for  $s\leq t\leq R$ and $R>0$ sufficiently large, we get

\begin{multline*}
J'(tv^{+}+sv^{-})(tv^{+})\leq m M(1)R^{2m}\xi_{1}(|\nabla (v^{+})|_{\Phi})^{2}\\
+mM(1) R^{2m} \xi_{1}(|\nabla (tv^{+})|_{\Phi})
\xi_{1}(|\nabla (v^{-})|_{\Phi}) \\
+mK_{1} R^{m}\xi_{1}(|\nabla (v^{+})|_{\Phi}) 
- \theta K_2 R^{\theta}|v^{+}|^{\theta}_{\theta}+ K_{3}\theta|\Omega|<0.
\end{multline*}

Arguing of the same way we get
$$
J'(tv^{+}+Rv^{-})(Rv^{-})<0, \ \ \mbox{for all} \ \ t\leq R.
$$

In particular,
$$
J'(rv^{+}+sv^{-})(rv^{+})>0 \ \ \mbox{and}  \ \
J'(tv^{+}+rv^{-})(rv^{-})>0, \ \ \mbox{for all} \ \ t,s \in [r,R]
$$
and
$$
J'(Rv^{+}+sv^{-})(Rv^{+})<0 \ \ \mbox{and}  \ \
J'(tv^{+}+Rv^{-})(Rv^{-})<0, \ \ \mbox{for all} \ \ t,s \in [r,R].
$$

Now the lemma follows applying Miranda's theorem \cite{Miranda}.
\hfill\rule{2mm}{2mm}

\hspace{0.5cm}

In the next Lemma we prove  monotonicity results for some functions that will be much useful in our arguments.

\begin{Le} \label{crescimentos}
\item[(a)] The function  \begin{eqnarray}\label{monotonicy1}
t\mapsto \widehat{M}(t)- \frac{1}{2}M(t)t \ \ \mbox{is
increasing}.
\end{eqnarray}

\item[(b)] For all $v \in W^{1,\Phi}_{0}(\Omega)$ with $v\geq 0$ and $v\neq 0$, the function 
\begin{eqnarray}\label{monotonicy1000}
t\mapsto \widehat{M}\biggl(\displaystyle\int_{\Omega}\Phi(tv)dx\biggl)- \frac{1}{2m}M\biggl(\displaystyle\int_{\Omega}\Phi(tv)dx\biggl)
\displaystyle\int_{\Omega}\phi(tv)(tv)^{2}dx \ \ \mbox{is
increasing}.
\end{eqnarray}

\item[(c)] The function 
\begin{eqnarray}\label{monotonicy2}
t\mapsto \frac{1}{2m}f(t)t- F(t) \ \ \mbox{is increasing, for all} \
\ |t|>0.
\end{eqnarray}
\item[(d)]
\begin{eqnarray}\label{dicaJoao}
\widehat{M}(t+s) \geq  \widehat{M}(t)+ \widehat{M}(s), \ \ \mbox{for all} \ \ t,s \in
[0,+\infty).
\end{eqnarray}
\end{Le}
\noindent {\bf Proof.} $(a)$ First of all,  let us
observe that, from $(M_{2})$ we have
\begin{eqnarray}\label{relacaoderivadafuncaoM}
M'(t)t\leq M(t), \ \ \mbox{for all} \ \ t\geq 0.
\end{eqnarray}
Now from (\ref{relacaoderivadafuncaoM}) we conclude $( \widehat{M}(t)- \frac{1}{2}M(t)t)'>0$ which implies
\begin{eqnarray*}
t\mapsto \widehat{M}(t)- \frac{1}{2}M(t)t \ \ \mbox{is
increasing}.
\end{eqnarray*}

Now we prove $(b)$. For each $v\in W^{1,\Phi}_{0}(\Omega)$ with $v\geq 0$ and $v\neq 0$, we define the function
\begin{eqnarray*}
\zeta(t)= \widehat{M}\biggl(\displaystyle\int_{\Omega}\Phi(tv)dx\biggl)- \frac{1}{2m}M\biggl(\displaystyle\int_{\Omega}\Phi(tv)dx\biggl)
\displaystyle\int_{\Omega}\phi(tv)(tv)^{2}dx.
\end{eqnarray*}
Then, recalling that $\Psi(t)=\phi(t)t$ ( see Lemma \ref{novomasseguefukagai} ), we have
\begin{multline*}
\zeta'(t)=M\biggl(\displaystyle\int_{\Omega}\Phi(tv)dx\biggl)\displaystyle\int_{\Omega}\Psi(tv)v dx\\
- \frac{1}{2m}M'\biggl(\displaystyle\int_{\Omega}\Phi(tv)dx\biggl)\displaystyle\int_{\Omega}\Psi(tv)v dx
\displaystyle\int_{\Omega}\Psi(tv)tv dx\\
-  \frac{1}{2m}M\biggl(\displaystyle\int_{\Omega}\Phi(tv)dx\biggl)\displaystyle\int_{\Omega}\bigl(\Psi'(tv)tv^{2}+\Psi(tv)v\bigl)dx.
\end{multline*}

Now using $(\phi_{3})$ we have
\begin{multline*}
\zeta'(t)\geq  \displaystyle\int_{\Omega}\Psi(tv)v dx
\biggl(M\biggl(\displaystyle\int_{\Omega}\Phi(tv)dx\biggl)\\
- \frac{1}{2m}M'\biggl(\displaystyle\int_{\Omega}\Phi(tv)dx\biggl)\displaystyle\int_{\Omega}\Psi(tv)tv dx\\
-  \frac{m-1}{2m}M\biggl(\displaystyle\int_{\Omega}\Phi(tv)dx\biggl)- \frac{1}{2m}
M\biggl(\displaystyle\int_{\Omega}\Phi(tv)dx\biggl)\biggl).
\end{multline*}

Using $(\phi_{2})$ we get

\begin{eqnarray*}
\zeta'(t)\geq && \displaystyle\int_{\Omega}\Psi(tv)v dx
\biggl(\frac{1}{2}M\biggl(\displaystyle\int_{\Omega}\Phi(tv)dx\biggl)-  \frac{1}{2}M'\biggl(\displaystyle\int_{\Omega}\Phi(tv)dx\biggl)
\displaystyle\int_{\Omega}\Phi(tv)dx\biggl).
\end{eqnarray*}
By (\ref{relacaoderivadafuncaoM}) we obtain $\zeta'(t)>0$ and the proof of $(b)$ is done.

Now we prove $(c)$. Using $(f_{4})$ we get
\begin{eqnarray}\label{relacaoderivadafuncaof}
f'(t)t\geq (2m-1)f(t),\ \ \mbox{for all} \ \ |t|\geq 0,
\end{eqnarray}
that implies
$$
(\frac{1}{2m}f(t)t- F(t))'>0, 
$$
that is,
\begin{eqnarray*}
t\mapsto \frac{1}{2m}f(t)t- F(t) \ \ \mbox{is increasing, for all} \
\ |t|>0.
\end{eqnarray*}

Finaly, we prove $(d)$. From $(M_{1})$ we obtain
\begin{eqnarray*}
\widehat{M}(t+s)&=& \displaystyle\int^{t+s}_{0}M(\tau) \ d\tau
=\widehat{M}(t)+ \displaystyle\int^{t+s}_{t}M(\tau) \ d\tau \nonumber\\
&=& \widehat{M}(t)+ \displaystyle\int^{s}_{0}M(\gamma+t) \ d\gamma
\nonumber
\\
&\geq & \widehat{M}(t)+ \displaystyle\int^{s}_{0}M(\gamma) \ d\gamma
\nonumber \\
&=& \widehat{M}(t)+ \widehat{M}(s), \ \ \mbox{for all} \ \ t,s \in
[0,+\infty).\hfill\rule{2mm}{2mm}
\end{eqnarray*}

\hspace{0.5cm}

Now, we can define a suitable  function and its  gradient vector
field which are related to functional $J$  and  will be involved in
particular in the application of the deformation lemma. Indeed, for
each $v\in W^{1,\Phi}_{0}(\Omega)$ with $v^{\pm}\neq 0$ we consider
$$
h^{v}:[0,+\infty)\times [0,+\infty)\rightarrow \mathbf{R} \ \
\mbox{given by}  \ \ h^{v}(t,s)=J(tv^{+}+sv^{-})
$$
and and its gradient $\Upsilon^{v}:[0,+\infty)\times
[0,+\infty)\rightarrow \mathbf{R}^{2}$ defined by
\begin{eqnarray*}
\Upsilon^{v}(t,s)&=&   \Big(\Upsilon^v_1(t,s),  \Upsilon^v_2(t,s) \Big) =
\Big(\frac{\partial h^{v}}{\partial t} (t,s),
\frac{\partial h^{v}}{\partial s} (t,s) \Big)\\
&=&(J'(tv^{+}+sv^{-})v^{+}, J'(tv^{+}+sv^{-})v^{-}),
\end{eqnarray*}
for every $(t,s) \in [0,+\infty)\times [0,+\infty)$. Furthermore, we
consider  the Hessian matrix of $h^v$ or  the Jacobian matrix of
$\Upsilon^v$, i.e.
\begin{displaymath}
 (\Upsilon^v)'(t,s)=
 \left( \begin{array}{ccc}
 \frac{\partial \Upsilon_1^v}{\partial t}(t,s) &   \frac{\partial \Upsilon_1^v}{\partial s}(t,s) \\
        \\
   \frac{\partial \Upsilon_2^v}{\partial t}(t,s) &  \frac{\partial \Upsilon_2^v}{\partial s}(t,s)
 \end{array} \right),
 \end{displaymath}
 for every $(t,s) \in [0,+\infty)\times [0,+\infty)$. Indeed,  in the following we aim to prove that,  if $w \in \mathcal{M}$, function $h^w$ has a critical point and in particular a  global maximum in $(t, s)= (1, 1)$,

\begin{Le} \label{resultadoscom11}
If $w\in \mathcal{M}$, then
\begin{description}
\item[(a)]
$$
h^{w}(t,s)<h^{w}(1,1)=J(w),
$$
for all $t,s\geq 0$ such that $(t,s)\neq (1,1)$.

\item[(b)]
$$
det(\Upsilon^{w})'(1,1)>0.
$$
\end{description}
\end{Le}
\noindent {\bf Proof.} Since $w\in \mathcal{M}$, then
$$
J'(w)w^{\pm}=J'(w^{+}+w^{-})w^{\pm}=0.
$$
Thus,
$$
\Upsilon^{w}(1,1)=(\frac{\partial h^{w}}{\partial t} (1,1),
\frac{\partial h^{w}}{\partial s} (1,1) )=(0,0).
$$

Moreover, from (\ref{estimativaporcima}) and (\ref{estimativaf2}), for $t$ and $s$ sufficiently large,  we
get
\begin{multline*}
h^{w}(t,s)=J(tw^{+}+sw^{-})\leq  m M(1)t^{2m}\xi_{1}(|\nabla (w^{+})|_{\Phi})^{2}+ m M(1)s^{2m}\xi_{1}(|\nabla (w^{-})|_{\Phi})^{2}\\
+mM(1) t^{2m} \xi_{1}(|\nabla (tw^{+})|_{\Phi})
\xi_{1}(|\nabla (w^{-})|_{\Phi}) \\+mM(1) s^{2m} \xi_{1}(|\nabla (sw^{-})|_{\Phi})
\xi_{1}(|\nabla (w^{+})|_{\Phi})\\
+ mK_{1} t^{m}\xi_{1}(|\nabla (w^{+})|_{\Phi}) + mK_{1} s^{m}\xi_{1}(|\nabla (w^{-})|_{\Phi})\\
-K_{2}\theta t^{\theta}|w^{+}|^{\theta}_{\theta}+ 2K_{3}\theta|\Omega|-K_{2}\theta s^{\theta}|w^{-}|^{\theta}_{\theta}.
\end{multline*}

Since $2m<\theta $, then
$$
\displaystyle\lim_{|(t,s)|\rightarrow +\infty}h^{w}(t,s)=-\infty,
$$
that implies $(1,1)$ is a critical point of $h^{w}$ and $h^{w}$ has
a global maximum point in $(a,b)$.

Now we prove that $a,b>0$. Suppose, by contradiction that $b=0$.
Thus, 
$$
J'(aw^{+})aw^{+}=0
$$ 
implies
\begin{eqnarray*}
M\biggl(\displaystyle\int_{\Omega} \Phi(a|\nabla w^{+}|)dx\biggl)\displaystyle\int_{\Omega} \phi(a|\nabla w^{+}|)a^{2}|\nabla w^{+}|^{2}dx=
\displaystyle\int_{\Omega}f(aw^{+})aw^{+} dx.
\end{eqnarray*}

Then, by Lemma \ref{novomasseguefukagai}, we have
\begin{eqnarray*}
M\biggl(\displaystyle\int_{\Omega} \Phi(a|\nabla w^{+}|)dx\biggl)\tau_{1}(a)a\displaystyle\int_{\Omega} 
\phi(|\nabla w^{+}|)|\nabla w^{+}|^{2}dx\geq
\displaystyle\int_{\Omega}f(aw^{+})aw^{+} dx, 
\end{eqnarray*}
which is equivalent
\begin{eqnarray*}
M\biggl(\displaystyle\int_{\Omega} \Phi(a|\nabla w^{+}|)dx\biggl)\displaystyle\int_{\Omega} \phi(|\nabla w^{+}|)|\nabla w^{+}|^{2}dx\geq
\displaystyle\int_{\Omega}\frac{f(aw^{+})aw^{+} dx}{\tau_{1}(a)a}. 
\end{eqnarray*}

Hence, 
\begin{eqnarray*}
\frac{M\biggl(\displaystyle\int_{\Omega} \Phi(a|\nabla w^{+}|)dx\biggl)}{\displaystyle\int_{\Omega} 
\Phi(a|\nabla w^{+}|)dx}\displaystyle\int_{\Omega} \Phi(a|\nabla w^{+}|)dx\displaystyle\int_{\Omega} \phi(|\nabla w^{+}|)
|\nabla w^{+}|^{2}dx\geq
\displaystyle\int_{\Omega}\frac{f(aw^{+})aw^{+} dx}{\tau_{1}(a)a}. 
\end{eqnarray*}

Since $\xi_{1}(a)=\tau_{1}(a)a$ we have
\begin{eqnarray}\label{comparar1}
\frac{M\biggl(\displaystyle\int_{\Omega} \Phi(a|\nabla w^{+}|)dx\biggl)}{\displaystyle\int_{\Omega} 
\Phi(a|\nabla w^{+}|)dx}\displaystyle\int_{\Omega} \Phi(|\nabla w^{+}|)dx
\displaystyle\int_{\Omega} \phi(|\nabla w^{+}|)|\nabla w^{+}|^{2}dx\geq
\displaystyle\int_{\Omega}\frac{f(aw^{+})aw^{+} dx}{\xi_{1}(a)^{2}}. 
\end{eqnarray}

On the other hand, since $J(w)w^{+}=0$ and $M$ is increasing, we get
\begin{eqnarray}\label{comparar2}
\frac{M\biggl(\displaystyle\int_{\Omega} \Phi(|\nabla w^{+}|)dx\biggl)}{\displaystyle\int_{\Omega} 
\Phi(|\nabla w^{+}|)dx}\displaystyle\int_{\Omega} \Phi(|\nabla w^{+}|)dx\displaystyle\int_{\Omega} 
\phi(|\nabla w^{+}|)|\nabla w^{+}|^{2}dx\leq
\displaystyle\int_{\Omega}f(w^{+})w^{+} dx. 
\end{eqnarray}

Considering (\ref{comparar1}) and (\ref{comparar2}) we have
\begin{eqnarray*}
&&\biggl[\frac{M\biggl(\displaystyle\int_{\Omega} \Phi(|\nabla w^{+}|)dx\biggl)}{\displaystyle\int_{\Omega} 
\Phi(|\nabla w^{+}|)dx}-\frac{M\biggl(\displaystyle\int_{\Omega} \Phi(a|\nabla w^{+}|)dx\biggl)}{\displaystyle\int_{\Omega} 
\Phi(a|\nabla w^{+}|)dx}\biggl]\displaystyle\int_{\Omega} 
\Phi(|\nabla w^{+}|)dx\displaystyle\int_{\Omega} \phi(|\nabla w^{+}|)|\nabla w^{+}|^{2}dx\\
&\leq &
\displaystyle\int_{\Omega}\biggl[\frac{f(w^{+})}{(w^{+})^{2m-1}}-\frac{f(aw^{+})a}{\xi_{1}(a)^{2}(w^{+})^{2m-1}}\biggl](w^{+})^{2m}dx.
\end{eqnarray*}

The last inequality,  $(M_{2})$ and $(f_{4})$ imply that $a\leq 1$, because otherwise we get
\begin{eqnarray*}
&&0<\biggl[\frac{M\biggl(\displaystyle\int_{\Omega} \Phi(|\nabla w^{+}|)dx\biggl)}{\displaystyle\int_{\Omega} 
\Phi(|\nabla w^{+}|)dx}-\frac{M\biggl(\displaystyle\int_{\Omega} \Phi(a|\nabla w^{+}|)dx\biggl)}{\displaystyle\int_{\Omega} 
\Phi(a|\nabla w^{+}|)dx}\biggl]\displaystyle\int_{\Omega} 
\Phi(|\nabla w^{+}|)dx\displaystyle\int_{\Omega} \phi(|\nabla w^{+}|)|\nabla w^{+}|^{2}dx\\
&\leq &
\displaystyle\int_{\Omega}\biggl[\frac{f(w^{+})}{(w^{+})^{2m-1}}-\frac{f(aw^{+})}{(aw^{+})^{2m-1}}\biggl](w^{+})^{2m}dx<0.
\end{eqnarray*}

Now note that
\begin{eqnarray}\label{aehmaximum}
h^{w}(a,0)&=&J(aw^{+})=J(aw^{+})-\frac{1}{2m}J'(aw^{+})(aw^{+})\nonumber\\
&=&\biggl[\widehat{M}\biggl(\displaystyle\int_{\Omega} \Phi(a|\nabla w^{+}|)dx\biggl)-
\frac{1}{2m}M\biggl(\displaystyle\int_{\Omega} \Phi(a|\nabla w^{+}|)dx\biggl)\displaystyle\int_{\Omega} 
\phi(a|\nabla w^{+}|)|aw^{+}|^{2}dx\biggl]
\nonumber\\
&+&\int_{\Omega}\biggl[\displaystyle\frac{1}{2m}f(aw^{+})aw^{+}-F(aw^{+})
 \biggl]dx.
\end{eqnarray}

Since $0<a\leq 1$ and using (\ref{monotonicy1000}) and  (\ref{monotonicy2}) in
(\ref{aehmaximum}) we obtain
\begin{eqnarray*}
h^{w}(a,0)&\leq &
\biggl[\widehat{M}\biggl(\displaystyle\int_{\Omega} \Phi(|\nabla w^{+}|)dx\biggl)-
\frac{1}{2m}M\biggl(\displaystyle\int_{\Omega} \Phi(|\nabla w^{+}|)dx\biggl)\displaystyle\int_{\Omega} \phi(|\nabla w^{+}|)
|w^{+}|^{2}dx\biggl]
\nonumber\\
&+&\int_{\Omega}\biggl[\displaystyle\frac{1}{2m}f(w^{+})w^{+}-F(w^{+})
\biggl]dx \\
&=& J(w^{+})-\frac{1}{2m}J'(w^{+})w= J(w^{+})=h^{w}(1,0).
\end{eqnarray*}

\noindent Now, our aim is to prove that
$$
J(w^+)=h^w(1,0) < J(w)=h^w(1,1).
$$

By Lemma \ref{limitacaoporbaixo1} we have that $J(w^{-})\geq 0$.
Thus,
\begin{eqnarray}\label{Sara}
J(w^+) \leq   J(w^{+}) + J(w^{-})=&&
 \widehat{M}\biggl(\displaystyle\int_{\Omega} \Phi(|\nabla w^{+}|)dx\biggl)
+\widehat{M}\biggl(\displaystyle\int_{\Omega} \Phi(|\nabla w^{-}|)dx\biggl)\nonumber\\
&-& \displaystyle\int_{\Omega} \left(F(w^{+})+F(w^{-})\right) dx.
\end{eqnarray}

By (\ref{dicaJoao}) we get
$$
J(w^+) <\widehat{M}\biggl(\displaystyle\int_{\Omega} \Phi(|\nabla w^{+}|)dx+ \displaystyle\int_{\Omega} \Phi(|\nabla w^{-}|)dx\biggl)-
\displaystyle\int_{\Omega} \left(F(w^{+})+F(w^{-}\right)) dx.
$$

Since the supports of $w^{+}$ and $w^{-}$ are disjoint, we obtain
$$
h^w(1,0)= J(w^+) <  J(w)=h^w(1,1),
$$
which is an absurd because $(a,0)$ is a maximum point. The same way
we prove that  $0<a$.\\

Now we will prove that $0< a,b\leq 1$. Since $(a,b)$ is another
critical point of $h^{w}$, we have
$$
M\biggl(\displaystyle\int_{\Omega} \Phi(a|\nabla  w^{+}|)dx+ \displaystyle\int_{\Omega} \Phi(b|\nabla w^{-}|)dx\biggl)
\displaystyle\int_{\Omega} \phi(a|\nabla  w^{+}|)|aw^{+}|^{2}dx=\displaystyle\int_{\Omega}f(aw^{+})aw^{+}
\ dx.
$$

Without loss of generality, we can suppose that $b\leq a$. Suppose by contradiction $a\geq1$. Thus
\begin{eqnarray*}
 \displaystyle\int_{\Omega}f(aw^{+})aw^{+}
\ dx \leq \frac{M\biggl(\displaystyle\int_{\Omega} \Phi(a|\nabla  w|)dx\biggl)}{\displaystyle\int_{\Omega} \Phi(a|\nabla  w|)dx}
\displaystyle\int_{\Omega} \Phi(a|\nabla  w|)dx\displaystyle\int_{\Omega} \phi(a|\nabla  w^{+}|)|aw^{+}|^{2}dx.
\end{eqnarray*}
Then,
\begin{eqnarray}\label{Giusepe}
 \displaystyle\int_{\Omega}\frac{f(aw^{+})}{(aw^{+})^{2m-1}}(w^{+})^{2m}
\ dx \leq \frac{M\biggl(\displaystyle\int_{\Omega} \Phi(a|\nabla  w|)dx\biggl)}{\displaystyle\int_{\Omega} \Phi(a|\nabla  w|)dx}
\displaystyle\int_{\Omega} \Phi(|\nabla  w|)dx\displaystyle\int_{\Omega} \phi(|\nabla  w^{+}|)|w^{+}|^{2}dx.
\end{eqnarray}

On the other hand, $J'(w)w^{+}=0$ implies
\begin{eqnarray}\label{Giusepe1}
\frac{M\biggl(\displaystyle\int_{\Omega} \Phi(|\nabla  w|)dx\biggl)}{\displaystyle\int_{\Omega} \Phi(|\nabla  w|)dx}
\displaystyle\int_{\Omega} \Phi(|\nabla  w|)dx
\displaystyle\int_{\Omega} \phi(|\nabla  w^{+}|)|w^{+}|^{2}dx=\displaystyle\int_{\Omega}f(w^{+})w^{+}
\ dx.
\end{eqnarray}

Combining (\ref{Giusepe}) and (\ref{Giusepe1}) and $(M_{2})$ and $(f_{4})$, we get
\begin{eqnarray}
&&0<\biggl[\frac{M\biggl(\displaystyle\int_{\Omega} \Phi(|\nabla  w|)dx\biggl)}{\displaystyle\int_{\Omega} \Phi(|\nabla  w|)dx}
-\frac{M\biggl(\displaystyle\int_{\Omega} \Phi(a|\nabla  w|)dx\biggl)}{\displaystyle\int_{\Omega} \Phi(a|\nabla  w|)dx}
\biggl]\displaystyle\int_{\Omega} \Phi(|\nabla  w|)dx
\displaystyle\int_{\Omega} \phi(|\nabla  w^{+}|)|w^{+}|^{2}dx \nonumber\\
&\leq &
\displaystyle\int_{\Omega}\biggl[\frac{f(w^{+})}{(w^{+})^{2m-1}}-\frac{f(aw^{+})}{(aw^{+})^{2m-1}}\biggl](w^{+})^{2m}
dx <0,
\end{eqnarray}

\noindent which is a contradiction. Then, $0<b\leq a < 1$.

Now we will prove that $h^{w}$ does not have global maximum in
$[0,1]\times [0,1]\backslash \{(1,1)\}$. We will show that
$$
h^{w}(a,b)< h^{w}(1,1).
$$

Note that $|\nabla (aw^{+}+bw^{-})|=a|\nabla w^{+}|+ b|\nabla w^{-}|\leq
|\nabla w^{+}|+|\nabla w^{-}|$ and since $\widehat{M}$ and $\Phi$ are increasing
we have
\begin{eqnarray}\label{***}
h^{w}(a,b)=J(aw^{+}+bw^{-})\leq
\widehat{M}\biggl(\displaystyle\int_{\Omega} \Phi(|\nabla  w|)dx\biggl)-\displaystyle\int_{\Omega}F(aw^{+}+bw^{-})
dx.
\end{eqnarray}

By $(f_{3})$ we get
$\displaystyle\frac{1}{2m}\displaystyle\int_{\Omega}f(aw^{+}+bw^{-})(aw^{+}+bw^{-})
dx\geq 0$. Thus, put this information in (\ref{***}) we obtain

\begin{eqnarray*}
h^{w}(a,b)&=& J(aw^{+}+bw^{-}) \leq
\widehat{M}\biggl(\displaystyle\int_{\Omega} \Phi(|\nabla  w|)dx\biggl)\\
&+&\displaystyle\int_{\Omega}\biggl[\frac{1}{2m}f(aw^{+}+bw^{-})(aw^{+}+bw^{-})-
F(aw^{+}+bw^{-})\biggl] dx.
\end{eqnarray*}

Since $w^{+}$ and $w^{-}$ have supports disjoint we get
\begin{eqnarray*}
h^{w}(a,b) \leq &&
\widehat{M}\biggl(\displaystyle\int_{\Omega} \Phi(|\nabla  w|)dx\biggl)\\
&+&\displaystyle\int_{\Omega}\biggl[\frac{1}{2m}f(aw^{+})(aw^{+})-
F(aw^{+})\biggl] dx\\
&+&  \displaystyle\int_{\Omega}\biggl[\frac{1}{2m}f(bw^{-})(bw^{-})-
F(bw^{-})\biggl] dx.
\end{eqnarray*}

Now, using  (\ref{monotonicy2}) we get 
\begin{eqnarray*}
&&\displaystyle\int_{\Omega}\biggl[\frac{1}{2m}f(aw^{+})(aw^{+})-
F(aw^{+})\biggl] dx\\
&+&  \displaystyle\int_{\Omega}\biggl[\frac{1}{2m}f(bw^{-})(bw^{-})-
F(bw^{-})\biggl] dx\\
&<& \displaystyle\int_{\Omega}\biggl[\frac{1}{2m}f(w^{+})(w^{+})-
F(w^{+})\biggl] dx\\
&+&  \displaystyle\int_{\Omega}\biggl[\frac{1}{2m}f(w^{-})(w^{-})-
F(w^{-})\biggl] dx
\end{eqnarray*}
which implies
$$
h^{w}(a, b)< J(w^{+}+w^{-})=J(w)=h^{w}(1,1)
$$
and item $(a)$ is proved.

Let us prove item $(b)$. Consider the notations
$\Upsilon^{w}_{1}(t,s)=J'(tw^{+}+ sw^{-})w^{+}$ and
$\Upsilon^{w}_{2}(t,s)=J'(tw^{+}+ sw^{-})w^{-}$. Thus,
\begin{eqnarray*}
\Upsilon^{w}_{1}(t,s) = M\biggl(\displaystyle\int_{\Omega}\Phi(t|\nabla w^{+}|+s|\nabla w^{-}|)dx\biggl)
\displaystyle\int_{\Omega}\phi(t|\nabla w^{+}|)t|\nabla w^{+}|^{2}dx-
\displaystyle\int_{\Omega}f(tw^{+})w^{+} dx
\end{eqnarray*}
and
\begin{eqnarray*}
\Upsilon^{w}_{2}(t,s)= M\biggl(\displaystyle\int_{\Omega}\Phi(t|\nabla w^{+}|+s|\nabla w^{-}|)dx\biggl)
\displaystyle\int_{\Omega}\phi(s|\nabla w^{-}|)s|\nabla w^{-}|^{2}dx -
\displaystyle\int_{\Omega}f(sw^{-})w^{-} dx.
\end{eqnarray*}
Then
\begin{eqnarray*}
&&\frac{\partial\Upsilon^{w}_{1}}{\partial
t}(t,s)\\
&=& M'\biggl(\displaystyle\int_{\Omega}\Phi(t|\nabla w^{+}|+s|\nabla w^{-}|)dx\biggl)
\displaystyle\int_{\Omega}\Psi(t|\nabla w^{+}|+s|\nabla w^{-}|)|\nabla w^{+}|dx
\displaystyle\int_{\Omega}\Psi(t|\nabla w^{+}|)|\nabla w^{+}|dx\\
&-&\displaystyle\int_{\Omega}f'(tw^{+})(w^{+})^{2}dx + M\biggl(\displaystyle\int_{\Omega}\Phi(t|\nabla w^{+}|+s|\nabla w^{-}|)dx\biggl)
\displaystyle\int_{\Omega}\Psi'(t|\nabla w^{+}|)|\nabla w^{+}|^{2}dx
\end{eqnarray*}
implies
\begin{eqnarray*}
\frac{\partial\Upsilon^{w}_{1}}{\partial
t}(1,1)\\
&=& M'\biggl(\displaystyle\int_{\Omega}\Phi(|\nabla w|)dx\biggl)
\displaystyle\int_{\Omega}\Psi(|\nabla w|)|\nabla w^{+}|dx
\displaystyle\int_{\Omega}\Psi(|\nabla w^{+}|)|\nabla w^{+}|dx\\
&-&\displaystyle\int_{\Omega}f'(w^{+})(w^{+})^{2}dx + M\biggl(\displaystyle\int_{\Omega}\Phi(|\nabla w|)dx\biggl)
\displaystyle\int_{\Omega}\Psi'(|\nabla w^{+}|)|\nabla w^{+}|^{2}dx\\
&\leq & m M'\biggl(\displaystyle\int_{\Omega}\Phi(|\nabla w|)dx\biggl)
\displaystyle\int_{\Omega}\Phi(|\nabla w|)dx
\displaystyle\int_{\Omega}\Psi(|\nabla w^{+}|)|\nabla w^{+}|dx\\
&-&\displaystyle\int_{\Omega}f'(w^{+})(w^{+})^{2}dx + M\biggl(\displaystyle\int_{\Omega}\Phi(|\nabla w|)dx\biggl)
\displaystyle\int_{\Omega}\Psi'(|\nabla w^{+}|)|\nabla w^{+}|^{2}dx.
\end{eqnarray*}

Using (\ref{relacaoderivadafuncaoM}) in the last equality we obtain
\begin{eqnarray*}
&&\frac{\partial\Upsilon^{w}_{1}}{\partial
t}(1,1)\leq m M\biggl(\displaystyle\int_{\Omega}\Phi(|\nabla w|)dx\biggl)
\displaystyle\int_{\Omega}\Psi(|\nabla w^{+}|)dx\\
&-&\displaystyle\int_{\Omega}f'(w^{+})(w^{+})^{2}dx + (m-1) M\biggl(\displaystyle\int_{\Omega}\Phi(|\nabla w|)dx\biggl)
\displaystyle\int_{\Omega}\Psi(|\nabla w^{+}|)dx\\
&=&(2m-1)\displaystyle\int_{\Omega}f(w^{+})w^{+}dx-\displaystyle\int_{\Omega}f'(w^{+})(w^{+})^{2}dx.
\end{eqnarray*}

From (\ref{relacaoderivadafuncaof}) we obtain
\begin{eqnarray}\label{determinante1}
\frac{\partial\Upsilon^{w}_{1}}{\partial t}(1,1)<0.
\end{eqnarray}

Arguing of the same way we conclude
\begin{eqnarray}\label{determinante2}
\frac{\partial\Upsilon^{w}_{2}}{\partial s}(1,1)<0.
\end{eqnarray}

Since $\frac{\partial\Upsilon^{w}_{1}}{\partial s}(1,1)=0$ and
$\frac{\partial\Upsilon^{w}_{2}}{\partial t}(1,1)=0$  and considering
(\ref{determinante1}) and (\ref{determinante2}) we conclude that
$det(\Upsilon^{w})'(1,1)=\frac{\partial\Upsilon^{w}_{1}}{\partial
t}(1,1)\frac{\partial\Upsilon^{w}_{2}}{\partial s}(1,1)>0$ and the item
$(b)$ is proved.\hfill\rule{2mm}{2mm}

\section{Proof of Theorem \ref{Theorem1}}

In this section we will prove the existence of $w \in \mathcal{M}$
in which the infimum of $J$ is attained on $ \mathcal{M}$. After,
following some arguments used in \cite{Alves1} by Alves and Souto (see also \cite{Weth}) and, in particular, applying a
deformation lemma,  we find that $w$ is a critical point of $J$ and
then a least energy nodal solution of $(P)$. In order to complete
the proof of Theorem \ref{Theorem1}, we conclude by showing that $w$
has exactly two nodal domains.

\medskip

\noindent First of all, by Lemma \ref{limitacaoporbaixo1}, there
exists $c_{0} \in \mathbf{R}$ such that
$$
0<c_{0}=\displaystyle\inf_{v \in \mathcal{M}}J(v).
$$

\noindent Thus, there exists a  minimizing sequence $(w_{n})$ in
$\mathcal{M}$ which  is bounded from Lemma \ref{limitacaoporbaixo1}.
again. Hence, by Sobolev Imbedding Theorem, without loss of
generality, we can assume up to a subsequence that there exist $w,
w_1, w_2 \in W^{1,\Phi}_{0}(\Omega)$ such that
\begin{eqnarray*}
w_{n}\rightharpoonup w, & w^+_{n}\rightharpoonup w_1,  &  w^-_{n}\rightharpoonup w_2  \quad  \mbox{in} \ \ W^{1,\Phi}_{0}(\Omega), \\
w_{n}\rightarrow w, &   w^+_{n}\rightarrow w_1,  &
w^-_{n}\rightarrow w_2    \quad \mbox{in $L^q(\Omega)$, $q \in (m,l^*)$}.
\end{eqnarray*}
Since  the transformations  $w \rightarrow w^+$ and $w \rightarrow
w^-$ are continuous from $L^q(\Omega)$ in $L^q(\Omega)$ (see Lemma
2.3 in \cite{Castro} with suitable adaptations),  we have that $w^+=
w_1 \geq 0$ and $w^- = w_2 \leq 0$. At this point, we can  prove
that $w \in \mathcal{M}$. Indeed, by  $w^+_{n}\rightarrow w^+$ and
$w^-_{n}\rightarrow w^-$ in $L^q(\Omega)$  it is, as $n \rightarrow
+ \infty$

$$
\displaystyle\int_{\Omega}|(w_{n})^{\pm}|^{q} dx \rightarrow
\displaystyle\int_{\Omega}|w^{\pm}|^{q} dx.
$$
Then, by Lemma \ref{limitacaoporbaixo2}, we conclude that
$w^{\pm}\neq 0$ and consequently $w= w^++w^-$ is sign-changing. By
Lemma  \ref{teovalormedio}, there exist  $t,s >0$ such that
\begin{eqnarray}\label{primeirots}
J'(tw^{+}+sw^{-})w^{+} =0, \nonumber  \\
J'(tw^{+}+sw^{-})w^{-} =0,
\end{eqnarray}
then $tw^{+}+sw^{-} \in \mathcal{M}$.  Now, let us prove that $t, s
\leq 1$. First let us observe that, since  $f$  has a quasicritical
growth, using compactness Lemma of Strauss \cite[Theorem A.I, p.338]{Strauss}, we obtain
$$
\displaystyle\int_{\Omega}f((w_{n})^{\pm})(w_{n})^{\pm} dx
\rightarrow \displaystyle\int_{\Omega}f(w^{\pm})w^{\pm} dx
$$
and
$$
\displaystyle\int_{\Omega}F((w_{n})^{\pm}) dx  \rightarrow
\displaystyle\int_{\Omega}F(w^{\pm}) dx.
$$

\medskip

\noindent Thus,  since $ J'(w_{n})w_{n}^{\pm} = 0$, by $(M_{1})$ we
have
\begin{eqnarray}\label{segundo}
J'(w^{+})w^{+}\leq 0 \ \ \mbox{and} \ \ J'(w^{-})w^{-}\leq 0.
\end{eqnarray}

\noindent Consequently, combining (\ref{primeirots}) and
(\ref{segundo}) and arguing as in the proof of Lemma
\ref{resultadoscom11} item $(a)$, we obtain $0<t,s \leq 1$.

\noindent In the next step we show that $J(tw^{+}+sw^{-})=c_{0}$ and
$t=s=1$ or better  $J(w)=c_{0}$. Indeed, since $t, s \leq 1 $ and
$w_n \rightharpoonup w$ as $n \rightarrow + \infty$,  exploiting the
arguments used in the proof of Lemma \ref{resultadoscom11} item
$(a)$ and the weak lower semicontinuity of both $J$ and $K$ with
$K(u)= J'(u) u$ on $W^{1,\Phi}_{0}(\Omega)$ described above we get
\begin{eqnarray*}
c_{0}\leq J(tw^{+}+sw^{-})&=& J(tw^{+}+sw^{-})-\frac{1}{2m}
J'(tw^{+}+sw^{-})(tw^{+}+sw^{-}) \\
&\leq & J(w^{+}+w^{-})-\frac{1}{2m}
J'(w^{+}+w^{-})(w^{+}+w^{-}) \\
&\leq & \displaystyle\liminf_{n\rightarrow +\infty}\biggl[
J(w^{+}_{n}+w^{-}_{n})-\frac{1}{2m}
J'(w^{+}_{n}+w^{-}_{n})(w^{+}_{n}+w^{-}_{n})\biggl] \\
&=& \displaystyle\lim_{n\rightarrow +\infty} J(w_n) =c_{0}.
\end{eqnarray*}

\noindent At this point, by using a quantitative deformation lemma
and adapting the arguments used  in  \cite{Weth} with slight
technical changes,  we point out  that $w$ is a critical point of
$J$, i.e. $J'(w)=0$ . If we reason by contradiction, we find that
there exist a positive constant $\alpha > 0$ and $v_0 \in
W^{1,\Phi}_{0}(\Omega)$, $\| v_0 \|=1$ such that
$$ J'(w) v_0  = 2 \alpha > 0. $$
By the continuity of $J'$, we can choose a radius $r > 0$ so that
\begin{equation*}
J'(v), v_0 =  \alpha > 0, \ \mbox{for every $v \in B_r(w) \subset
W^{1,\Phi}_{0}(\Omega)$ with $v^{\pm} \neq 0$.}
\end{equation*}
Let us fix $D= (\xi, \chi) \times (\xi, \chi) \subset \mathbf{R}^2$
with $0 < \xi < 1 < \chi$ such that
\begin{itemize}
\item[(i)] $(1, 1) \in  D$ and $\Upsilon^{w}(t,s)=(0, 0) $ in $\overline{D}$ if and only if $(t, s)=(1, 1)$;

\item[(ii)] $c_0 \notin h^{w}(\partial D)$;

\item[(iii)] $\left\{ t w_0^+ + s w_0^- : (t, s ) \in \overline{D}   \right\} \subset B_r(w)$,
\end{itemize}
where $h^{w}$ and $\Upsilon^{w}$ are defined as  in Section 3  and
satisfy Lemma \ref{resultadoscom11}. At this point, we can choose a
smaller  radius $r' > 0$ such that
\begin{equation}\label{r'}
\mathcal{B}= \overline{B_{r'}(w)} \subset B_r(w) \ \mbox{and
$\mathcal{B} \cap \left\{ t w^+ + s w^- : (t, s ) \in \partial D
\right\} = \emptyset$.}
\end{equation}
Now define a continuous mapping $\rho: W^{1,\Phi}_{0}(\Omega) \rightarrow
[0, + \infty)$ such that
$$\rho(u):= \operatorname{dist}(u, \mathcal{B}^c), \ \mbox{for all $u \in W^{1,\Phi}_{0}(\Omega)$},  $$
then a bounded Lipschitz vector field $V: W^{1,\Phi}_{0}(\Omega)
\rightarrow W^{1,\Phi}_{0}(\Omega)$ given by
$$ V(u)= - \rho(u) v_0 $$
and, for every $u \in W^{1,\Phi}_{0}(\Omega)$, denoting  by $\eta(\tau)=
\eta(\tau, u)$ we consider  the following Cauchy problem
\begin{eqnarray*}
\begin{cases}
\eta'(\tau)= V(\eta(\tau)), \ \mbox{for all $\tau > 0$,} \\
\eta(0)=u.
\end{cases}
\end{eqnarray*}
Now, we observe that there exist a continuous deformation
$\eta(\tau, u)$ and $\tau_0 > 0$ such that for all $\tau \in [0,
\tau_0]$ the following properties hold:

\begin{itemize}
\item[(a)] $\eta(\tau, u)=u$  for all $u \notin \mathcal{B}$;
\item[(b)]  $\tau \rightarrow J(\eta(\tau, u))$ is decreasing  for all $\eta(\tau, u) \in \mathcal{B}$;
\item[(c)] $J(\eta(\tau, w_0)) \leq J(w) - \frac{r' \alpha}{2} \tau$.
\end{itemize}

\noindent Item $(a)$ follows immediately by the definition of $\rho$. Indeed,  $u \notin \mathcal{B}$ implies  
$\rho(u)=0$ and the unique solution satisfying the above Cauchy problem is constant with constant value $u$. As concerns as item $(b)$, let us first observe that, since
$\eta(\tau) \in \mathcal{B} \subset B_r(w)$,  $ J'(\eta(\tau))  v_0
= \alpha > 0$ and, by definition of $\rho$, it is $\rho(\eta(\tau))
> 0$. Now, differentiating  $ J$ with respect to $\tau$, for all
$\eta(\tau) \in \mathcal{B}$, we have that
\begin{eqnarray*}
\frac{d}{d \tau} \left(J(\eta(\tau) \right) =   J'(\eta(\tau))
\eta'(\tau)  =   - \rho(\eta(\tau)) J'(\eta(\tau))  v_0  =
-\rho(\eta(\tau)) \alpha < 0
\end{eqnarray*}
thus concluding that $J(\eta(\tau, u))$ is decreasing with respect to $\tau$. \\
In order to prove item $(c)$, being $\tau_0 > 0$ such that
$\eta(\tau, u) \in \mathcal B$  for every $ 0 \leq \tau \leq
\tau_0$,  we can assume without loss of generality
$$\| \eta(\tau, w)-w  \| \leq \frac{r'}{2}   \Longleftrightarrow  \eta(\tau, w) \in \overline{B_{\frac{r'}{2}}(w)}, \ \mbox{for every $ 0 \leq \tau \leq \tau_0$}. $$
Thus, since $\rho(\eta(\tau, w))= \operatorname{dist}( \eta(\tau,
w), \mathcal{B}^c) \geq \frac{r'}{2}$ it follows that
$$ \frac{d}{d \tau} \left(J(\eta(\tau, w) \right) = -\rho(\eta(\tau, w)) \alpha \leq \frac{r' \alpha}{2}   $$
and, integrating in $[0, \tau_0]$ we finally get
$$ J(\eta(\tau, w_0)) -  J(w) \leq - \frac{r' \alpha}{2} \tau. $$

\medskip

\noindent At this point, let us consider  a suitable deformed path
$\overline{\eta}_0: \overline{D} \rightarrow X$ defined by
$$\overline{\eta}_{\tau_0}(t, s):= \eta(\tau_0, tw^+ +s w^-), \ \mbox{ for all $(t, s) \in \overline D$}   $$
so that
$$ \max_{(t, s) \in \overline D} J(\overline{\eta}_{\tau_0}(t, s)) < c_0. $$
Indeed,  by $(b)$ and the fact that $\eta$ satisfies the initial
condition $\eta(0, u)= u$,  for all $(t,s) \in \overline D - \{(1,
1)\}$ it is
\begin{eqnarray*}
J(\overline{\eta}_{\tau_0}(t, s)) &=&  J(\eta(\tau_0, tw^+ +s w^-)) \leq J(\eta(0, tw^+ +s w^-)) \\
& =&  J(tw^+ +s w^-)= h^w(t, s) < c_0,
\end{eqnarray*}
and, for $(t, s)= (1, 1)$, by $(c)$ we get
 \begin{eqnarray*}
 J(\overline{\eta}_{\tau_0}(1, 1))& =&  J(\eta(\tau_0, w^+ + w^-))=  J(\eta(\tau_0, w)) \\
& \leq & J(w) - \frac{r' \alpha}{2} \tau_0 < J(w) < c_0.
\end{eqnarray*}
Then, $\overline{\eta}_{\tau_0}(\overline{D}) \cap \mathcal{M} \neq
\emptyset $, i.e.
\begin{equation}\label{notM}
\overline{\eta}_{\tau_0}(t,s) \notin \mathcal{M}, \ \mbox{for all
$(t,s) \in \overline D$.}
\end{equation}
On the other side, defined $\Lambda_{\tau_0}: \overline D \rightarrow
\mathbf{R}^2$ such that
$$ \Lambda_{\tau_0}:= \Big(\frac{ J'(\overline{\eta}_{\tau_0}(t, s)) (\overline{\eta}_{\tau_0}(t, s))^+ }{t},
\frac{ J'(\overline{\eta}_{\tau_0}(t, s))
(\overline{\eta}_{\tau_0}(t, s))^- }{s}     \Big), $$ we observe
that, for all $(t, s) \in \partial D$, by \eqref{r'} and $(a)$  for
$\tau= \tau_0$,  it is
$$ \Lambda_{\tau_0}(t, s)= \Big( J'(t w^+ +s w^-) w^+ ,   J'(t w^+ +s w^-), w^-    \Big) = \Phi^w(t, s).  $$
Then, since  by   Brouwer's topological degree
$$\deg(\Lambda_{\tau_0}, D, (0, 0))=  \deg(\Phi^w, D, (0, 0))= \operatorname{sgn}(\det(\Phi^w)'(1, 1))=1, $$
we get that $\Lambda_{\tau_0}$ has a zero $(\overline t, \overline s)
\in D$ namely
$$\Lambda_{\tau_0}(\overline t, \overline s)=(0, 0) \Longleftrightarrow   J'(\overline{\eta}_{\tau_0}(\overline t, \overline s))
 (\overline{\eta}_{\tau_0}(\overline t, \overline s))^{\pm} =0. $$
Consequently there exists $(\overline t, \overline s) \in  D $ such
that $\overline{\eta}_{\tau_0}(\overline t, \overline s) \in
\mathcal{M}$ and we have a contradiction with \eqref{notM}. We
conclude that $w$ is a critical point of $J$.

\bigskip

\noindent Finally, we prove that $w$ has exactly two nodal domains
or equivalently it changes sign exactly once. Let us observe that
assumptions  $(M_1)$, $(f_1)$ and $(f_2)$ ensure  that  $w$ is
continuous and then $\widetilde{\Omega}= \{x \in \Omega: w(x) \neq 0
\}$ is open.  Suppose by contradiction that $\widetilde{\Omega}$ has
more than two components or $w$ has more than two nodal domains and,
since  $w$ changes sign, without loss of generality, we can assume 
$$w= w_1+w_2+w_3, \ \mbox{where $w_1 \geq 0$, $w_2 \leq 0$, $w_3 \neq 0$,} $$
and
$$\operatorname{supp}(w_i) \cap \operatorname{supp}(w_j) = \emptyset, \ \mbox{for $i \neq j$, $i, j=1,2,3$.} $$
So the disjointness of the supports  combined with  $J'(w)=0$
implies
$$ \left\langle J'(w_1+w_2), w_1 \right\rangle = 0 = \left\langle J'(w_1+w_2), w_2 \right\rangle .    $$
Since $ 0 \neq w_1= (w_1+w_2)^+$ and $ 0 \neq w_2= (w_1+w_2)^-$, by previous arguments,  there exist $t,s \in (0,1]$ such that $t (w_1+w_2)^+ + s (w_1+w_2)^+ \in \mathcal{M} $ namely  $t w_1 + s w_2 \in \mathcal{M}$ and then $J(t w_1 + s w_2) \geq c_0$. \\
On the other side,  $0 \neq w_3 \in \mathcal{N}$, Lemma 2.1 (i) and
the  arguments used  in the proof of Lemma \ref{resultadoscom11}
imply that
$$J(t w_1 + s w_2) \leq J( w_1 +  w_2) <   J( w_1 +  w_2) + J(w_3)= J(w) = c_0$$
then a contradiction and we conclude that $w_3 =0$.  Thus,  the
proof of Theorem \ref{Theorem1}      is complete.\hfill\rule{2mm}{2mm}


\begin{thebibliography}{99}

\bibitem{adams} A. Adams and J. F. Fournier, Sobolev spaces, {\it 2nd ed.}, Academic Press, (2003).





\bibitem{alvescorrea} C.O. Alves and  F.J.S.A. Corr\^{e}a , {\it On
existence of solutions for a class of problem involving a nonlinear
operator}, Comm. Appl. Nonlinear Anal., 8(2001)43-56.

\bibitem{alvescorreama}
C.O. Alves, F.J.S.A. Corr\^{e}a  and T.F. Ma,  {\it Positive
solutions for a quasilinear elliptic equation of Kirchhoff type}
{\it Comput. Math. Appl., 49(2005)85-93}.

\bibitem{AlvesFigueiredo}
C. O. Alves  and G. M. Figueiredo  {\it Nonlinear perturbations of a
periodic Kirchhoff equation in $\mathbf{R}^{N}$.} {\it Nonlinear
Anal., 75(2012)2750-2759}.

 \bibitem{Alves1} C. O. Alves and M. A S. Souto, {\it Existence of least energy nodal solution for a
Schr\"{o}dinger-Poisson system in bounded domains}, Zeitschrift fur Angewandte Mathematik und Physik (Printed ed.), p. 1153-1166, 2013.

\bibitem{Arosio} A. Arosio, {\it On the nonlinear Timoshenko-Kirchoff beam equation},
Chin. Ann. Math., 20 (1999), 495-506.

\bibitem{Arosio1} A. Arosio, {\it A geometrical nonlinear correction to the Timoshenko beam equation},
 Nonlinear Anal. 47(2001), 729-740.








\bibitem{Azzollini} A. Azzollini, {\it The elliptic Kirchhoff equation in $\mathbf{R}^{N}$ perturbed by a local nonlinearity.}
{\it Differential Integral Equations 25 (2012), no. 5-6, 543-554.}

\bibitem{Weth} T. Bartsch, T. Weth and  M. Willem, {\it Partial symmetry of least energy
nodal solution to some variational problems}, J. Anal. Math. 1,
(2005)1-18.

\bibitem{Bartsch} T. Bartsch and T. Weth, {\it Three nodal solutions of singularly perturbed
elliptic equations on domains without topology}, Ann. Inst. H.
Poincar\'e Anal. Non Lin\'eaire 22 (2005), 259-281.

\bibitem{Bartsch1} T. Bartsch, Z. Liu and T. Weth, {\it Sign changing solutions of superlinear
Schr\"{o}dinger equations}, Comm. Partial Differential Equations 29
(2004), 25-42.

\bibitem{Strauss} H. Berestycki and P.L. Lions, {\it  Nonlinear scalar field equations, I - existence
of a ground state}, Arch. Rat. Mech. Analysis, 82, (1983), 313-346.

  

\bibitem{Brezis_Lieb} H. Brezis. and E. Lieb, {\it A relation between pointwise convergence of functions and
convergence of functinals}, Proc. Amer. Math. Soc. 88 (1983),
486-490.

\bibitem{Cammaroto} F. Cammaroto  and L. Vilasi, {\it On a Schr\"{o}dinger-Kirchhoff-type equation involving the
$p(x)$-Laplacian}, Nonlinear Anal. 81 (2013) 42-53.

\bibitem{Castro} A. Castro, J. Cossio and J. Neuberger, {\it A sign-changing solution for a
superlinear Dirichlet problem}, Rocky Mountain J. Math. 27, 4
(1997), 1041-1053.


\bibitem{Cheng} B. Cheng, X. Wu and J. Liu, {\it Multiple solutions for a class of Kirchhoff
type problems with concave nonlinearity}, NoDEA Nonlinear
Differential Equations Appl. 19 (2012), 521-537.

\bibitem{Chen} C. Chen, H. Song and  Z. Xiu, {\it Multiple solutions for p-Kirchhoff equations in
$\mathbf{R}^{N}$}, Nonlinear Anal. 86 (2013) 146-156.


\bibitem{SChen} S. Chen and L. Li, {\it Multiple solutions for the nonhomogeneous Kirchhoff equation
on $\mathbf{R}^{N}$}, Nonlinear Anal. Real World Appl., 14 (2013)
1477-1486.


\bibitem{Chung} N. T. Chung, {\it Existence of solutions for a class of Kirchhoff type problems in Orlicz-Sobolev spaces}, 
Acta Univers. Apulensis, 37(2014), 111-123. 

\bibitem{Chung1} N. T. Chung, {\it Multiple solutions for a nonlocal problem in Orlicz-Sobolev spaces}, 
 Ric. Mat. 63 (2014), no. 1, 169-182.

\bibitem{Chung2} N. T. Chung, {\it Three solutions for a class of nonlocal problems in Orlicz-Sobolev spaces}, 
 J. Korean Math. 50(2013), 1257-1269. 


\bibitem{Donaldson2}{T.K. Donaldson and N.S. Trudinger,}{\it Orlicz-Sobolev spaces and imbedding theorems,} 
J. Funct. Anal. 8 (1971) 52-75.


\bibitem{jmaa} G.M. Figueiredo, {\it Existence of positive solution for a Kirchhoff  problem type with
critical growth via truncation argument}, {\it J. Math. Anal. Appl.
401 (2013), 706-713.}


\bibitem{Cristian} G.M. Figueiredo, C. Morales, J. R. Santos Junior and A. Siarez, {\it 
Study of a nonlinear Kirchhoff equation with non-homogeneous material }, {\it J. Math. Anal. Appl.
416(2014), 597-608.}

\bibitem{Giovany1} G.M. Figueiredo and J. R. Santos Junior, {\it 
Multiplicity and concentration behavior of positive solutions for a Schrodinger-Kirchhoff type problem 
via penalization method }, {\it ESAIM: Control, Otimiz. Calc. variat. vol 20, issue 2, (2014), 389-415.}

\bibitem{Giovany2} G.M. Figueiredo, N. Ikoma and J. R. Santos Junior, {\it 
Existence and concentration result for the Kirchhoff type equations with general nonlinearities }, 
{\it ARMA, 213 issue 3(2014), 931-979}.
 
\bibitem{Giovany3} G.M. Figueiredo R. G. Nascimento, {\it 
Existence of a nodal solution with minimal energy for a Kirchhoff equation }, 
{\it Math. Nachrichten, 288 (2015), pg 48-60.}

\bibitem{Giovany4} G.M. Figueiredo J. A. Santos, {\it 
On a nonlocal multivalued problem  in an Orlicz-Sobolev space via Krasnoselskii's genus }, 
{\it J. Convex Anal. 22 ( 2015)-447-446.}

\bibitem{Giovany5} G.M. Figueiredo J. A. Santos, {\it 
On a Phi-Kirchhoff multivalued problem with critical growth in an Orlicz-Sobolev space}, 
{\it Asympt. Anal. 89(2014), 151-172 . }
 
 

\bibitem{fukagai}N. Fukagai and K. Narukawa,
{\it Positive solutons of quasilinear elliptic equations with
critical Orlicz-Sobolev nonlinearity on $\mathbf{R}^{N}$},
Funkciallaj Ekvacioj, 49(2006)235-267.

\bibitem{Gossez} J.P. Gossez, {\it Orlicz-Sobolev spaces and nonlinear elliptic boundary
value problems.} In: Fuk, Svatopluk and Kufner, Alois (eds.):
Nonlinear Analysis, Function Spaces and Applications, Proceedings of
a Spring School held in Horn Bradlo, 1978. [Vol 1]. BSB B. G.
Teubner Verlagsgesellschaft, Leipzig, 1979. pp. 59-94.

\bibitem{He}
X. He  and W. Zou,  {\it Existence and concentration of positive
solutions for a Kirchhoff equation in $\mathbf{R}^{3}$}. {\it J.
Differential Equations, 252(2012)1813-1834}.



\bibitem{kirchhoff}
G. Kirchhoff, {\it Mechanik, Teubner,Leipzig, 1883}.


\bibitem{Li}
Y. Li, F. Li  and J. Shi,  {\it Existence of a positive solution to
Kirchhoff type problems without compactness conditions}.  {\it J.
Differential Equations, 253(2012)2285-2294}.

\bibitem{Liao}
YJia-Feng Liao, Peng Zhang, Jiu Liu and Chun-Lei Tang {\it Existence and multiplicity of positive solutions for a class of Kirchhoff 
type problems with singularity}.  {\it J. Math. Anal. Appl. 430 (2015), 1124-1148}
  

\bibitem{Liu} X. Liu and Y. Sun, {\it Multiple positive solutions for Kirchhoff type problems with
singularity}, Commun. Pure Appl. Anal. 12 (2013), 721-733.

\bibitem{Liu1} Z. Liu and S. Guo, {\it Existence of positive ground state solutions for Kirchhoff type problems.}, 
Nonlinear Anal. 120 (2015), 1-13.



\bibitem{Liang} S. Liang and  Shaoyun Shi, {\it Soliton solutions to Kirchhoff type problems involving the
critical growth in $\mathbf{R}^{N}$}, Nonlinear Anal. 81 (2013)
31-41.

\bibitem{ma}
T.F. Ma, {\it Remarks on an elliptic equation of Kirchhoff type}.
{\it Nonlinear Anal.,  63, 5-7(2005)1967-1977}.

\bibitem{Mao} A. Mao and Z. Zhang, {\it Sign-changing and multiple solutions of Kirchhoff type problems
without the P.S. condition}, Nonlinear Anal. 70 (2009) 1275-1287.

\bibitem{Mao1} A. Mao and Shixia Luan,  {\it Sign-changing solutions of a class of nonlocal quasilinear elliptic
boundary value problems}, J. Math. Anal. Appl. 383 (2011) 239-243.

\bibitem{Miranda}
C. Miranda, {\it Un' osservazione su un teorema di Brouwer}, Boll.
Un. Mat. Ital., 3 (1940) 5-7.

\bibitem{Naimen}
D. Naimen, {\it On the Brezis-Nirenberg problem with a Kirchhoff type perturbation.},Adv. Nonlinear Stud. 15 (2015), 135-156.

\bibitem{Shuai}
W. Shuai, {\it Sign-changing solutions for a class of Kirchhoff-type problem in bounded domains}, J. Diff. Equations, 254 (2015) 1256-1274. 


\bibitem{Sun} J. Sun, {\it Resonance problems for Kirchhoff type
equations}, Discrete Contin. Dyn. Syst., 33 issue 5(2013), 2139-2154.

\bibitem{Xiang} M. Xiang, B. Zhang and X. Guo, {\it Infinitely many solutions for a fractional Kirchhoff type 
problem via Fountain Theorem}, Nonlinear Anal. 120 (2015),  299-313.


\bibitem{Wang1}
J. Wang , L. Tian  , J. Xu  and F. Zhang,  {\it Multiplicity and
concentration of positive solutions for a Kirchhoff type problem
with critical growth}. {\it J. Differential Equations,
253(2012)2314-2351}.

\bibitem{Wang2} L. Wang,   {\it On a quasilinear Schr\"{o}dinger-Kirchhoff-type equation with radial potentials.} 
{\it Nonlinear Anal. 83 (2013),
58-68.}



\bibitem{Wu}
Xian Wu, {\it Existence of nontrivial solutions and high energy
solutions for Schr\"{o}dinger-Kirchhoff-type equations in
$\mathbf{R}^{N}$}. {\it Nonlinear Anal. Real World Appl.,
12(2011)1278-1287.}

\bibitem{Zhang} Z. Zhang and K. Perera, {\it Sign changing solutions of Kirchhoff type problems
via invariant sets of descent flow}, J. Math. Anal. Appl. 317 (2006)
456-463.


\end{thebibliography}
\end{document}